\documentclass[11pt]{amsart}
\input xy
\xyoption{all}
\usepackage{amsmath,amsthm, amscd, amssymb, amsfonts}

\usepackage[usenames, dvipsnames]{xcolor}

\newcommand{\comp}{\mbox{comp}}
\newcommand{\coev}{\mbox{coev}}
\newcommand{\ev}{\mbox{ev}}

\newcommand{\rhob}{\bar{\rho}}

\newcommand{\otb}{{\overline{\otimes}}}
\newcommand{\otk}{{\otimes_{\ku}}}
\newcommand{\Mo}{{\mathcal M}}
\newcommand{\No}{{\mathcal N}}

\newcommand{\Bimo}{{\mathcal Bimod}}

\newcommand{\ca}{{\mathcal C}}

\newcommand{\ot}{{\otimes}}

\newcommand{\kc}{{\mathcal K}}

\newcommand{\Sc}{{\mathcal S}}
\newcommand{\Zc}{{\mathcal Z}}

\newcommand{\cha}{{\mathcal A}}

\newcommand{\Do}{{\mathcal D}}

\newcommand{\Bc}{{\mathcal B}}

\newcommand{\Fc}{{\mathcal F}}

\newcommand{\YD}{{\mathcal YD}}

\newcommand{\op}{\rm{op}}

\newcommand{\ra}{\rm{ra}}
\newcommand{\la}{\rm{la}}
\newcommand{\cf}{\rm{CF}}

\newcommand{\ku}{{\Bbbk}}

\newcommand{\uno}{ \mathbf{1}}
\newcommand{\C}{{\mathcal C}}

\newcommand{\id}{\mbox{\rm id\,}}
\newcommand{\fp}{\mbox{\rm FPdim\,}}
\newcommand{\Id}{\mbox{\rm Id\,}}

\newcommand{\Ab}{\mbox{\rm Ab\,}}

\newcommand{\Bimod}{\mbox{\rm Bimod\,}}

\newcommand{\vect}{\mbox{\rm vect\,}}

\newcommand{\Rex}{\mbox{\rm Rex\,}}
\newcommand{\Fun}{\operatorname{Fun}}
\newcommand\Rep{\operatorname{Rep}}

\newcommand\Hom{\operatorname{Hom}}
\newcommand\uhom{\underline{\Hom}}

\newcommand{\End}{\operatorname{End}}

\renewcommand{\_}[1]{\mbox{$_{\left( #1 \right)}$}}

\theoremstyle{plain}

\numberwithin{equation}{section}

\newtheorem{teo}{Theorem}[section]

\newtheorem{lema}[teo]{Lemma}

\newtheorem{cor}[teo]{Corollary}

\newtheorem{prop}[teo]{Proposition}

\theoremstyle{definition}

\newtheorem{defi}[teo]{Definition}

  \newtheorem{exa}[teo]{Example}

\theoremstyle{remark}

\newtheorem{rmk}[teo]{Remark}

\def\pf{\begin{proof}}

\def\epf{\end{proof}}

\theoremstyle{remark}

\begin{document}

\title[The Character algebra for Hopf algebras]
{The Character algebra for module categories over Hopf algebras}
\author[   Bortolussi and Mombelli  ]{ Noelia Bortolussi and Mart\'in Mombelli
 }

\address{Facultad de Matem\'atica, Astronom\'\i a y F\'\i sica
\newline \indent
Universidad Nacional de C\'ordoba
\newline
\indent CIEM -- CONICET
\newline \indent Medina Allende s/n
\newline
\indent (5000) Ciudad Universitaria, C\'ordoba, Argentina}
 \email{ bortolussinb@gmail.com, nbortolussi@famaf.unc.edu.ar }
\email{martin10090@gmail.com, mombelli@mate.uncor.edu
\newline \indent\emph{URL:}\/ http://www.mate.uncor.edu/$\sim$mombelli/welcome.html}

\begin{abstract} Given a finite dimensional Hopf algebra $H$ and an exact indecomposable module category $\Mo$ over $\Rep(H)$, we explicitly compute the adjoint algebra $\cha_\Mo$ as an object in the category of Yetter-Drinfeld modules over $H$, and the space of class functions $\cf(\Mo)$ associated to $\Mo$, as introduced by K. Shimizu \cite{Sh2}. We  use our construction to describe these algebras when $H$ is a group algebra and a dual group algebra. This result allows us to compute the adjoint algebra for certain  group-theoretical fusion categories.
\end{abstract}

\date{\today}
\maketitle


\section*{Introduction}

In the paper \cite{Sh1},  the author introduces the notion of adjoint algebra $\cha_\ca$ and the space of class functions $\cf(\ca)$ for an arbitrary finite tensor category $\ca$. The adjoint algebra is defined as the end $  \int_{X\in \ca} X\ot X^*$. The dual object $\cha^*_\ca$ is a crucial ingredient in Lyubashenko's theory of the modular group action in non semisimple tensor categories \cite{L1}, \cite{L2}.

Both, the adjoint algebra and the space of class functions, are  interesting objects that generalize the well known adjoint representation and the character algebra of a finite group. In \cite{Sh1} many results concerning table of characters, conjugacy classes, and orthogonality relations of characters in finite group theory have been generalized to the setting of fusion categories.  Also, in \cite{Sh3} the adjoint algebra was used to develop a theory of integrals for finite tensor categories.

Assume $\Mo$ is an arbitrary module categories over a finite tensor category $\ca$.   In \cite{Sh2}, the author introduces the notion of adjoint algebra $\cha_\Mo$ and the space of class functions $\cf(\Mo)$ associated to $\Mo$, generalizing the definitions given in \cite{Sh1}. The main task of  this paper is the explicit computation of those objects in the particular case $\ca$ is the representation category of a finite dimensional Hopf algebra.

\medbreak

Assume that $\ca$ is a finite tensor category, and $\Mo$ an exact  left $\ca$-module category with action functor $\otb:\ca\times \Mo\to \Mo$. Then, we can consider  the  functor $\rho_\Mo:\ca\to \Rex(\Mo)$, $\rho_\Mo(X)(M)=X\otb M$, $X\in \ca$, $M\in\Mo$. Here $\Rex(\Mo)$ denotes the category of right exact endofunctors of $\Mo$. The right adjoint of the action functor is $\rho_\Mo^{\ra}:  \Rex(\Mo)\to\ca$, explicitly described as
$$ \rho_\Mo^{\ra}(F)= \int_{M\in \Mo} \uhom(M, F(M)),$$
for any  $F\in \Rex(\Mo)$ \cite[Thm. 3.4]{Sh2}.  Here for any $M\in \Mo$, $\uhom(M, -) $ is the right adjoint of the functor $\ca\to \Mo$, $X\mapsto X\otb M$. It is called the \textit{internal Hom }of the module category $\Mo$. The adjoint algebra is defined as $\cha_\Mo= \rho_\Mo^{\ra}(\Id_\Mo)$. This object has a half-braiding $\sigma_{\Mo}(X):  \cha_{\Mo} \otimes X \to X\ot  \cha_\Mo$ defined as the unique morphism in $\ca$ such that the diagram

\begin{equation*}
  \xymatrix@C=34pt@R=34pt{
    \cha_\Mo \otimes X
    \ar[rr]^-{\pi_{\Mo}(X \otb M) \otimes \id_V}
    \ar[d]_{\sigma_\Mo(X)}
    & & \underline{\End}(X \otb M) \otimes X
    \ar[d]^{\simeq} \\
    X \otimes  \cha_\Mo 
    \ar[r]^-{\id_X \otimes \pi_\Mo}
    & X \otimes \uhom(M, M)
    \ar[r]^{\simeq }
    & \uhom(X \otb M, M)
  }
\end{equation*}
is commutative. Here $\pi_{\Mo}(M):\cha_\Mo \to   \uhom(M, M)$ is the dinatural transformation of the end $\cha_\Mo$. Turns out that $(\cha_\Mo, \sigma_\Mo)$ is a commutative algebra in the Drinfeld center $\Zc(\ca)$. Although  this description of the half-braiding of $\cha_\Mo$  is rather  clear, for us it was complicated to use it to make calculations in particular examples. However, there is another way of describing this structure.

If $\Bc$ is a $\ca$-bimodule category, one can consider the \textit{relative center }$\Zc_\ca(\Bc)$. When $\ca$ is considered as a bimodule over itself, the relative center coincides with the Drinfeld center. The correspondence $\Bc \mapsto\Zc_\ca(\Bc) $  is in fact part of a 2-functor
$$ \Zc_\ca: {}_\ca\Bimod \to \Ab_\ku,$$
where ${}_\ca\Bimod $ is the 2-category of finite $\ca$-bimodule categories, bimodule functors and bimodule natural transformations, and $\Ab_\ku$ is the 2-category of finite abelian $\ku$-linear categories.
Both categories, $\Rex(\Mo)$ and $\ca$ are $\ca$-bimodule categories. Turns out that $ \rho_\Mo^{\ra}$ has a $\ca$-bimodule structure \cite[Section 3.4]{Sh2}. Applying the 2-functor $\Zc_\ca$ one obtains a functor $\Zc_\ca(\rho_\Mo^{\ra}): \Zc_\ca(\Rex(\Mo))\simeq \ca^*_\Mo\to \Zc(\ca)$. Hence $(\cha_\Mo,\sigma_\Mo)=\Zc_\ca(\rho_\Mo^{\ra})(\Id_\Mo)$.
\medbreak

Assume $H$ is a finite dimensional Hopf algebra. If $\Mo$ is an exact indecomposable module category over $\Rep(H)$, we describe explicitly the adjoint algebra $(\cha_\Mo,\sigma_\Mo)$ and the space of class functions $\cf(\Mo)$. For this purpose, we need to explain all ingredients in the construction of those objects. Our description of both algebras relies heavily on the explicit description of module categories over Hopf algebras. In Section \ref{subsection:modcat-over-hopf} we embark on this task. Module categories over $\Rep(H)$  are categories ${}_K\Mo$ of finite dimensional left $K$-modules, where $K$ are certain $H$-comodule algebras. We also recall how to describe module functor categories, and that there is a monoidal equivalence $\Rep(H)^*_{{}_K\Mo}\simeq {}_K^H\Mo_K$. This equivalence will be used when explaining the functor $\Zc_\ca(\rho_\Mo^{\ra})$. Another ingredient is the internal Hom. In this section we also describe, in a precise way, the internal Hom of the module category ${}_K\Mo$.  In Section \ref{Section:character-algebramain}, after recalling the definitions of \cite{Sh2}, for an object  $P\in {}_K^H\Mo_K$, representing a module functor in $F_P\in \End_{\Rep(H)}({}_K\Mo)$, we explictly give the structure of the functor
$$ F_P\mapsto \int_{M\in \Mo} \uhom(M, F_P( M)).$$
For this we compute, in an explicit way, the end $\int_{M\in \Mo} \uhom(M, F_P( M))$ as an object in the category ${}^H_H\YD$ of Yetter-Drinfeld modules over $H$. 
In Section \ref{Section:example-comp}, we illustrate this description in the particular cases when $H$ is a group algebra or  its dual. As a direct consequence, we compute the adjoint algebra and the space of class functions for certain group-theoretical fusion categories.

\subsection*{Acknowledgments} This  work  was  partially supported by CONICET and Secyt (UNC),  Argentina. We also thank the referee for her/his comments that improve the presentation of the paper.

\section{Preliminaries}

Let $\ku$ be an algebraically closed field.
All algebras are assumed to be over $\ku$. If $A$ is an algebra, we shall denote by ${}_A\Mo$ (respectively $\Mo_A$) the category of finite dimensional left $A$-modules (respectively right $A$-modules). If $A, B$ are two algebras, we shall denote by ${}_B\Mo_A$ the category of finite dimensional $(B,A)$-modules.
From now on, all categories are assumed to be abelian $\ku$-linear, and all functors are $\ku$-linear.
\subsection{Hopf algebras}\label{Section:prelim-hopf} For a Hopf algebra $H$, we shall denote by $\Delta:H\to H\otk H$ the comultiplication, $\Sc:H\to H$ the antipode, and $\epsilon:H\to \ku$ the counit. We shall use Sweedler's notation: $\Delta(h)=h\_1\ot h\_2$, $h\in H$. The category ${}_H\Mo$ has a canonical structure of tensor category with monoidal product given by $\otk$. We shall denote this tensor category by $\Rep(H).$

For a finite dimensional Hopf algebra $H$, we shall denote by ${}^H_H\YD$ the category of finite-dimensional \textit{Yetter-Drinfeld modules.}  An object $V\in {}^H_H\YD$ is a left $H$-module $\cdot:H\otk V\to V$, and a left $H$-comodule $\lambda:V\to H\otk V$ such that
\begin{equation}\label{yd-module} \lambda(h\cdot v)=h\_1 v\_{-1} \Sc(h\_3)\ot h\_2\cdot v\_0,
\end{equation}
for any $h\in H, v\in V$.  If $V\in {}^H_H\YD$, the map
$\sigma_X: V\otk X\to X\otk V$, given by $\sigma_X(v\ot x)=v\_{-1}\cdot x\ot v\_0$ is a half-braiding for $V$, and this correspondence establishes a monoidal equivalence ${}^H_H\YD\simeq \Zc(\Rep(H))$.

\subsection{Finite categories}
A category $\ca$ is
\emph{finite} \cite{EO} if
\begin{itemize}
\item it has finitely many simple objects;
 \item each simple
object $X$ has a projective cover $P(X)$; \item the $\Hom$ spaces
are finite-dimensional; \item each object has finite length.
\end{itemize}
Equivalently, a category is finite if it is equivalent to a category ${}_A\Mo$ for some finite dimensional algebra $A$.

If $\Mo, \No$ are two finite categories, and $F:\Mo\to\No$ is a functor, we shall denote by $F^{\la}, F^{\ra}: \No\to \Mo$, its left adjoint, respectively right adjoint of $F$, if it exists. We shall also denote by $\Rex(\Mo,\No)$ the category of right exact functors from $\Mo$ to $\No$.
\subsection{Ends and coends}

We briefly recall the notion of end and coend. The reader is referred to \cite{Mac}. Let $\ca$, $\Do$ be categories, and let $S, T: \ca^{\op} \times \ca \to \Do$ be functors. A {\em dinatural transformation} $\xi: S \xrightarrow{..} T$ is a collection of morphisms in $\Do$
\begin{align*}
  \xi_{X}: S(X,X) \to T(X,X), \quad X \in \ca,
\end{align*}
such that for any morphism $f: X \to Y$ in $\ca$
\begin{equation}\label{dinaturality1}
  T(\id_X, f) \circ \xi_X \circ S(f, \id_X) = T(f, \id_Y) \circ \xi_Y \circ S(\id_Y, f).
\end{equation}
An {\em end} of $S$ is a pair $(E, p)$ consisting of an object $E \in \Do$ and a dinatural transformation $p: E \xrightarrow{..} S$ satisfying the following universal property.
For any pair $(D, q)$ consisting of an object $D\in \Do$ and a dinatural transformation $q: D \xrightarrow{..} S$, there exists a \textit{unique} morphism $h: D\to E$ in $\Do$ such that $q_X = p_X \circ h$ for any $X \in \ca$. A {\em coend} of $S$ is the dual notion of an end, this means that it is a pair $(C, \pi)$ consisting of an object $C \in \Do$ and a dinatural transformation $\pi: S \xrightarrow{..} C$ with the following universal property. For any pair $(B, t)$, where $B\in \Do$ is an object and $t: S\xrightarrow{..} B$ is a dinatural transformation, there exists a unique morphism $h:C\to B$ such that $h\circ \pi_X=t_X$ for any $X \in \ca$.

The end  and  coend of the functor $S$ are denoted, respectively, as
$$
  \int_{X \in \mathcal{C}} S(X,X)
  \quad \text{and} \quad
  \int^{X \in \mathcal{C}} S(X,X). $$

\section{Representations of tensor categories}\label{Section:represent}

For basic notions on finite tensor categories we refer to \cite{EGNO}, \cite{EO}. Let $\C$ be a finite tensor category over $\ku$. A (left) \emph{module} over 
$\ca$ is a  finite  category $\Mo$ together with a $\ku$-bilinear 
bifunctor $\otb: \ca \times \Mo \to \Mo$, exact in each variable,  endowed with 
 natural associativity
and unit isomorphisms 
$$m_{X,Y,M}: (X\otimes Y)\otb M \to X \otb
(Y\otb M), \ \ \ell_M: \uno \otb M\to M.$$ These isomorphisms are subject to the following conditions:
\begin{equation}\label{left-modulecat1} m_{X, Y, Z\otb M}\; m_{X\otimes Y, Z,
M}= (\id_{X}\otb m_{Y,Z, M})\;  m_{X, Y\otimes Z, M}(\alpha_{X,Y,Z}\otb\id_M),
\end{equation}
\begin{equation}\label{left-modulecat2} (\id_{X}\otb l_M)m_{X,{\bf
1} ,M}= r_X \otb \id_M,
\end{equation} for any $X, Y, Z\in\C, M\in\Mo.$ Here $\alpha$ is the associativity constraint of $\C$.
Sometimes we shall also say  that $\Mo$ is a $\ca$-\emph{module} or a $\ca$-\emph{module category}.

\medbreak

Let $\Mo$ and $\Mo'$ be a pair of $\C$-modules. A\emph{ module functor} is a pair $(F,c)$, where  $F:\Mo\to\Mo'$  is a functor equipped with natural isomorphisms
$$c_{X,M}: F(X\otb M)\to
X\otb F(M),$$ $X\in  \ca$, $M\in \Mo$,  such that
for any $X, Y\in
\ca$, $M\in \Mo$:
\begin{align}\label{modfunctor1}
(\id_X \otb  c_{Y,M})c_{X,Y\otb M}F(m_{X,Y,M}) &=
m_{X,Y,F(M)}\, c_{X\otimes Y,M}
\\\label{modfunctor2}
\ell_{F(M)} \,c_{\uno ,M} &=F(\ell_{M}).
\end{align}

There is a composition
of module functors: if $\Mo''$ is another $\C$-module and
$(G,d): \Mo' \to \Mo''$ is another module functor then the
composition
\begin{equation}\label{modfunctor-comp}
(G\circ F, e): \Mo \to \Mo'', \qquad  e_{X,M} = d_{X,F(M)}\circ
G(c_{X,M}),
\end{equation} is
also a module functor.

\smallbreak  

A \textit{natural module transformation} between  module functors $(F,c)$ and $(G,d)$ is a 
 natural transformation $\theta: F \to G$ such
that for any $X\in \ca$, $M\in \Mo$:
\begin{gather}
\label{modfunctor3} d_{X,M}\theta_{X\otb M} =
(\id_{X}\otb \theta_{M})c_{X,M}.
\end{gather}
 Two module functors $F, G$ are \emph{equivalent} if there exists a natural module isomorphism
$\theta:F \to G$. We denote by $\Fun_{\ca}(\Mo, \Mo')$ the category whose
objects are module functors $(F, c)$ from $\Mo$ to $\Mo'$ and arrows module natural transformations. 

\smallbreak
Two $\C$-modules $\Mo$ and $\Mo'$ are {\em equivalent} if there exist module functors $F:\Mo\to
\Mo'$, $G:\Mo'\to \Mo$, and natural module isomorphisms
$\Id_{\Mo'} \to F\circ G$, $\Id_{\Mo} \to G\circ F$.

A module is
{\em indecomposable} if it is not equivalent to a direct sum of
two non trivial modules. Recall from \cite{EO}, that  a
module $\Mo$ is \emph{exact} if   for any
projective object
$P\in \ca$ the object $P\otb M$ is projective in $\Mo$, for all
$M\in\Mo$. If $\Mo$ is an exact indecomposable module category over $\ca$, the dual category $\ca^*_\Mo=\End_\ca(\Mo)$ is a finite tensor category \cite{EO}. The tensor product is the composition of module functors.

 A \emph{right module category} over $\ca$
 is a finite  category $\Mo$ equipped with an exact
bifunctor $\otb:  \Mo\times  \ca\to \Mo$ and natural   isomorphisms 
$$\widetilde{m}_{M, X,Y}: M\otb (X\ot Y)\to (M\otb X) \otb Y, \quad r_M:M\otb \uno\to M$$ such that
\begin{equation}\label{right-modulecat1} \widetilde{m}_{M\otb X, Y ,Z }\; \widetilde{m}_{M,X ,Y\ot Z } (\id_M \otb a_{X,Y,Z})=
(\widetilde{m}_{M,X , Y}\ot \id_Z)\, \widetilde{m}_{M,X\ot Y ,Z },
\end{equation}
\begin{equation}\label{right-modulecat2} (r_M\otb \id_X)  \widetilde{m}_{M,\uno, X}= \id_M\otb l_X.
\end{equation}

If $\Mo,  \Mo'$ are right $\ca$-modules, a module functor from $\Mo$ to $  \Mo'$ is a pair $(T, d)$ where
$T:\Mo \to \Mo'$ is a  functor and $d_{M,X}:T(M\otb X)\to T(M)\otb X$ are natural  isomorphisms
such that for any $X, Y\in
\ca$, $M\in \Mo$:
\begin{align}\label{modfunctor11}
( d_{M,X}\ot \id_Y)d_{M\otb X, Y}T(m_{M, X, Y}) &=
m_{T(M), X,Y}\, d_{M, X\ot Y},
\\\label{modfunctor21}
r_{T(M)} \,d_{ M,\uno} &=T(r_{M}).
\end{align}

Assume that $\Mo,  \No$ are  categories, $F:\Mo\to \No$ is a  functor with right adjoint $G:\No\to\Mo$. We shall denote by $\epsilon:F\circ G\to \Id_{\No}$, $\eta:\Id_{\Mo}\to G\circ F$, the counit and unit of the adjunction. The next result will be needed later.
\begin{lema}\label{modfunct-adjoint} \cite[ Lemma 2.11]{DSS} The following holds.
\begin{itemize}
\item[(i)] If $\Mo,  \No$ are left $\ca$-module categories  and $(F,c):\Mo\to \No$ is a module functor then $G$ has a module functor structure given by
$$ e^{-1}_{X,N}=G(\id_X\otb \epsilon_N) G(c_{X,G(N)}) \eta_{X\otb G(N)},$$
for any $X\in \ca$, $N\in \No$.
\item[(ii)]   If $\Mo,  \No$ are right $\ca$-module categories  and $(F,d):\Mo\to \No$ is a module functor then $G$ has a module functor structure given by
$$ h^{-1}_{N,X}=G(\epsilon_N\otb \id_X) G(d_{G(N),X}) \eta_{G(N)\otb X},$$
for any $X\in \ca$, $N\in \No$.\qed
\end{itemize}
\end{lema}

\subsection{Bimodule categories} Assume $\Do$ is another finite tensor category. A $(\ca, \Do)-$\emph{bimodule category}  is a category $\Mo$  with left $\ca$-module category
 and right $\Do$-module category  structure
together with natural
isomorphisms 
\begin{align}\label{associativ-constraint-bimodulecat}
\gamma_{X,M,Y}:(X\otb M) \otb Y\to X\otb (M\otb Y),
\end{align} 
$X\in\ca, Y\in\Do, M\in \Mo$, satisfying 
\begin{align}\label{bimod11}
\xymatrix{
((X\ot Y)\otb M)\otb Z \ar[d]_{m^l\otimes \id} \ar[r]^{\gamma} & (X\ot Y)\otb(M\otb Z) \ar[dd]^{m^l} \\ 
(X\otb(Y\otb M))\otb Z \ar[d]_{\gamma} \\ 
X\otb((Y\otb M)\otb Z) \ar[r]_{\id\otimes \gamma} & X\otb (Y\otb (M\otb Z)), }
\;
\end{align} 

\begin{align}\label{bimod12}
\xymatrix{
(X\otb M)\otb (Y\ot Z) \ar[d]_{m^r} \ar[r]^{\gamma} & X\otb (M\otb (Y\ot Z)) \ar[dd]^{id\otimes m^r} \\ 
((X\otb M)\otb Y)\ot Z \ar[d]_{\gamma\otimes \id} \\ 
(X\otb (M\otb Y))\ot Z \ar[r]_{\gamma} & X\otb ((M\otb Y)\otb Z),
},
\end{align} 
\begin{align}\label{bimod13}
\xymatrix{
(\uno \otb M)\otb \uno\ar[r]^{\gamma_{\uno, M, \uno}}\ar[d]_{l_M \ot \id_\uno} & 1\otb (M\otb \uno)\ar[dd]^{r_{M}}\\
M\otb \uno \ar[d]_{\id_\uno\ot r_M}&\\
M & \uno\otb M \ar[l]^{l_M},}
\end{align} 
where $m^l$ and $m^r$ are the associativity isomorphisms of the left, respectively right, module category. If $\Mo, \No$ are $(\ca, \Do)-$bimodule categories, a \textit{bimodule functor} is a triple $(F,c,d):\Mo\to \No$, where $(F,c)$ is a $\ca$-module functor, $(F,d)$ is a $\Do$-module functor and equation
\begin{equation}\label{bimod-funct}
\gamma_{X, F(M),Y}(c_{X,M}\otb\id_Y) d_{X\otb M, Y}=(\id_X\otb d_{M,Y})c_{X, M\otb Y} F(\gamma_{X,M,Y}),
\end{equation}
is fulfilled for any $X\in \ca$, $Y\in \Do$, $M\in \Mo$.

It is known that   $\Mo$ is a $(\ca, \Do)-$bimodule category if and only if it is a left  $\ca\boxtimes\Do^{\op}$-module category, and a bimodule functor is the same as a $\ca\boxtimes\Do^{\op}$-module functor. See for example \cite{Gr2}.

If $\Mo, \No$ are left $\C$-module categories, then $\Rex(\Mo, \No)$ is a $\C$-bimodule category as follows. If $X\in \ca$, $F\in \Rex(\Mo, \No)$, $M\in \Mo$, then 
\begin{equation}\label{bimod-action-funct} (X\otb F)(M)=X\otb F(M), \quad (F\otb X)(M)=F(X\otb M).
\end{equation}

\subsection{The internal Hom}\label{subsection:internal hom} Let $\ca$ be a  tensor category and $\Mo$ be  a left $\C$-module category. For any pair of objects $M, N\in\Mo$, the \emph{internal Hom} is an object $\uhom(M,N)\in \C$ representing the functor $\Hom_{\Mo}(-\otb M,N):\ca\to \vect_\ku$. This means that there are natural isomorphisms, one the inverse of each other,
\begin{equation}\label{Hom-interno}\begin{split}\phi^X_{M,N}:\Hom_{\ca}(X,\uhom(M,N))\to \Hom_{\Mo}(X\otb M,N), \\
\psi^X_{M,N}:\Hom_{\Mo}(X\otb M,N)\to \Hom_{\ca}(X,\uhom(M,N)),
\end{split}
\end{equation}
 for all $M, N\in \Mo$, $X\in\ca$.
Sometimes we shall denote the internal Hom of the module category $\Mo$ by $\uhom_\Mo$ to emphasize that it is related to this module category.

For any $X\in \ca$, $M, N\in \Mo$ define 
$$\coev^\Mo_{X, M}:X\to  \uhom(M,X\otb M), \quad \ev^\Mo_{M, N}:\uhom(M,N)\otb M\to N,$$
$$\coev^\Mo_{X, M}=\psi^X_{M,X\otb M}(\id_{X\otb M}),  \quad\ev^\Mo_{M, N}= \phi^{\uhom(M,N)}_{M,N}(\id_{\uhom(M,N)}).$$

Define also $f_M=\ev^\Mo_{M,M}(\id_{\uhom(M,M)}\otb \ev^\Mo_{M, M})$, and
$$\comp^\Mo_M:\uhom(M,M)\ot \uhom(M,M)\to \uhom(M,M),$$
$$\comp^\Mo_M=\psi^{\uhom(M,M)\ot \uhom(M,M)}_{M,M}(f_M). $$
It is known, see \cite{EO}, that $\uhom(M,M)$ is an algebra in the category $\ca$ with product given by $\comp^\Mo_M.$
\subsection{The relative center}\label{Section:relative center}

Let $\ca$ be a tensor category and $\Mo$ a $\ca$-bimodule category. The \emph{
relative center} of $\Mo$ is the category of $\ca$-bimodule functors from $\ca$ to $\Mo$. We denote the relative center of $\Mo$ by $\Zc_\ca(\Mo)$ .
Explicitly, objects of  $\Zc_\ca(\Mo)$ are pairs $(M,\sigma)$, where $M$ is an object of $\Mo$ and
$$\sigma_X:M\overline{\otimes}X\xrightarrow{\sim}X\overline{\otimes} M $$ is a family of natural   isomorphisms such that
\begin{equation}\label{half-braid}
m^l_{X,Y,M}\sigma_{X\otimes Y}=(\id_X\otb \sigma_Y)\gamma_{X,M,Y} (\sigma_X\otb \id_Y)  m^r_{M,X,Y},
\end{equation}
where $\gamma_{X,M,Y}:(X\overline{\otimes} M)\overline{\otimes} Y \to X\overline{\otimes} (M\overline{\otimes} Y)$ are the associativity constraints of the left and right actions on $\Mo$, see 
\eqref{associativ-constraint-bimodulecat}. The isomorphism $\sigma$ is called the \textit{half-braiding} for $M$.

As  explained in \cite[Section 3.6]{Sh2}, the relative center can be thought of as a 2-functor
$$\Zc_\ca: {}_\ca\Bimo \to \Ab_\ku,$$
where ${}_\ca\Bimo$ is the 2-category whose 0-cells are $\ca$-bimodule categories, 1-cells are bimodule functors and 2-cells are bimodule natural transformations. Also $\Ab_\ku$ is the 2-category of finite $\ku$-linear abelian categories. If $\Mo, \No$ are $\ca$-bimodule categories, then $\Zc_\ca(\Mo)$ is the relative center. If $(F,c,d):\Mo\to \No$ is a bimodule functor, then $\Zc_\ca(F):\Zc_\ca(\Mo)\to \Zc_\ca(\No)$ is the functor $\Zc_\ca(F)(M,\sigma)=(F(M), \widetilde\sigma)$, where $\widetilde\sigma_X: F(M)\otb X\to X\otb F(M)$ is defined as
\begin{equation}\label{relative-bimod-funct} 
\widetilde\sigma_X=c_{X,M} F(\sigma_X)d^{-1}_{M,X},
\end{equation}
for any $X\in \ca$. 

The following example is \cite[Example 3.11]{Sh2}.
\begin{exa}\label{relat-center-rex} If $\Mo, \No$ are  exact $\ca$-module categories, then $\Rex(\Mo,\No)$ is a $\ca$-bimodule category, see \eqref{bimod-action-funct}. In this case there exists an equivalence $\Zc_\ca(\Rex(\Mo,\No))\simeq \Fun_\ca(\Mo,\No)$.
\end{exa}

\begin{exa} When $\ca$ is considered as a $\ca$-bimodule category, then $\Zc_\ca(\ca)=\Zc(\ca)$ is the usual center of the category $\ca$.
\end{exa}

\begin{rmk}\label{center:module functor} If $(X,\sigma)\in \Zc(\ca)$ and $\Mo$ is a left $\ca$-module category, then the functor  $L_{(X,\sigma)}:\Mo\to \Mo$ given by $L_{(X,\sigma)}(M)=X\otb M$ is a $\ca$-module functor. The module structure is given by
$$c^{(X,\sigma)}_{Y,M}: X\otb (Y\otb M)\to Y\otb (X\otb M),$$
$$c^{(X,\sigma)}_{Y,M}=m_{Y,X,M}(\sigma_Y\otb \id_M) m^{-1}_{X,Y,M}, $$
for any $X, Y\in \ca$, $M\in \Mo$.
\end{rmk}

\begin{defi}\label{forget-center} For any exact indecomposable left $\ca$-module category $\Mo$, we shall denote by 
$$\Fc_\Mo:\Zc(\ca^*_\Mo)\to \ca^*_\Mo,\quad (V,\sigma)\mapsto V,$$
the forgetful functor. In particular $\Fc_\ca: \Zc_\ca(\ca)\to \ca$ is the usual forgetful functor.
\end{defi}

\subsection{Morita invariance of the Drinfeld center}\label{Section:moritainvariance}
Let $\Mo$ be an exact indecomposable module category over $\ca$. Using results of P. Schauenburg \cite{Sch}, K. Shimizu proved in \cite[Section 3.7]{Sh2} that there exists a braided monoidal equivalence
$$\theta_\Mo:\Zc(\ca)\to \Zc(\ca^*_\Mo).$$
For later uses, we shall recall the definition of this equivalence. Let $(V,\sigma)\in \Zc(\ca)$. Then $\theta_\Mo(V,\sigma):\Mo\to \Mo$ is the functor defined as $\theta_\Mo(V,\sigma)(M)=V\otb M$, for any $M\in\Mo.$ The module structure of the functor $\theta_\Mo(V,\sigma)$ is 
$$c_{X,M}:\theta_\Mo(V,\sigma)(X\otb M)=V\otb (X\otb M) \to X\otb (V\otb M),$$
given by the composition
$$V\otb (X\otb M)\xrightarrow{m^{-1}_{V,X,M}} (V\ot X)\otb M \xrightarrow{\sigma_X\otb \id}   (X\ot V)\otb M \xrightarrow{m_{X,V,M}}   X\otb (V\otb M).$$
Then $\theta_\Mo(V,\sigma)$ becomes a $\ca$-module functor. It remains to explain how the functor $\theta_\Mo(V,\sigma)$ is an object in the center of  $\ca^*_\Mo$. For any $(F,d)\in \ca^*_\Mo$ we have to define a half-braiding $\tau_{(F,d)}:\theta_\Mo(V,\sigma)\circ (F,d)\to (F,d)\circ\theta_\Mo(V,\sigma).$ This is the module natural transformation defined by
\begin{equation}\label{braiding-center}
(\tau_{(F,d)})_M:V\otb F(M)\to F(V\otb M), \quad 
(\tau_{(F,d)})_M=d^{-1}_{V,M},
\end{equation}
for any $M\in \Mo$.

\section{Module categories over Hopf algebras}\label{subsection:modcat-over-hopf} Throughout this section $H$ will denote a finite dimensional Hopf algebra. We shall present families of module categories over $\Rep(H)$, and compute explicitly its internal Hom and their module functor categories.

\medbreak

If $\lambda:K\to H\otk K$ is a left $H$-comodule algebra then the
category of finite-dimensional left $K$-modules ${}_K\Mo$ is a
module category over $\Rep(H)$ with action $\otb:\Rep(H)\times
{}_K\Mo\to {}_K\Mo$, $X\otb M=X\otk M$, for all $X\in \Rep(H),
M\in {}_K\Mo$. The left $K$-module structure on $X\otk M$ is given
by $\lambda$, that is, if $k\in K$, $x\in X, m\in M$ then 
$$k\cdot
(x\ot m)= \lambda(k) (x\ot m)=k\_{-1}\cdot x\ot k\_{0}\cdot m.$$

\begin{teo}\cite[Prop.1.20]{AM}\label{sub-coid}
If $K$ is right $H$-simple then ${}_K\Mo$ is an exact indecomposable module category over $\Rep(H)$. Moreover, if $\Mo$ is an exact indecomposable module category over $\Rep(H)$, there exists a right $H$-simple left $H$-comodule algebra $K$, with trivial coinvariants, such that $\Mo\simeq {}_K\Mo$.\qed
\end{teo}

\begin{rmk} If $K, S$ are isomorphic $H$-comodule algebras, then the categories ${}_K\Mo, {}_S\Mo$ are equivalent as $\Rep(H)$-module categories. The converse is not always true.
\end{rmk}
\subsection{The internal Hom}\label{Section:homint-hf} We shall explicitly compute the internal Hom of module categories over $\Rep(H)$.

If $M, N$ are left $K$-modules, then the space $\Hom_K(H\otk M, N)$ has a left $H$-action given by
$$(h\cdot \alpha)(t\ot m)=\alpha(th\ot m),$$
for any $h,t\in H$, $\alpha\in \Hom_K(H\otk M, N)$, $m\in M$. 
We can identify the space $\Hom_K(H\otk M, N)$ with the subspace of $H^*\otk \Hom_\ku(M,N)$ consisting of elements $\sum_i f_i\ot T_i\in H^*\otk \Hom_\ku(M,N)$ such that
\begin{equation}\label{internal-id}
\sum_i  <f_i, k\_{-1}h>  T_i(k\_0\cdot m)=\sum_i <f_i, h>k\cdot T_i(m),
\end{equation}
for any $h\in H$, $k\in K$, $m\in M$. An element $\sum_i f_i\ot T_i $  is seen as a map from $H\otk M$ to $N$, sending $h\ot m$ to $\sum_i <f_i,h>  T_i(m).$ We shall freely use this identification from now on. Condition \eqref{internal-id} says that this morphism is a $K$-module map.

For any $K$-module $M$, the space $\Hom_K(H\otk M, M)$ has an algebra structure as follows. If $\sum_i f_i\ot T_i,$ $\sum_j g_j\ot U_j $ are elements in  $\Hom_K(H\otk M, M)$, the product is defined by
\begin{equation}\label{product-internal-hom}
\big(\sum_i f_i\ot T_i\big)\big(\sum_j g_j\ot U_j \big)=\sum_{i,j}  f_ig_j\ot T_i\circ  U_j.
\end{equation}
The proof of the next result is straightforward.
\begin{lema}\label{product-stab} With the product described in \eqref{product-internal-hom}, $\Hom_K(H\otk M, M)$ becomes an $H$-module algebra.\qed
\end{lema}

\begin{lema}\label{int-hom-hopf} Let $M, N\in{}_K\Mo$, and $\uhom(M,N)$ the  internal Hom of the module category  ${}_K\Mo$. There is an isomorphism of $H$-modules  
$$\uhom(M,N)\simeq \Hom_K(H\otk M, N).$$
When $M=N$ this isomorphism is an $H$-module algebra isomoprhism.
\end{lema}
\pf Let $X\in \Rep(H)$. The maps 
$$\phi:\Hom_H(X, \Hom_K(H\otk M, N))\to \Hom_K(X\otk M, N),$$
$$\psi:\Hom_K(X\otk M, N)\to \Hom_H(X, \Hom_K(H\otk M, N)),$$
defined by
$\phi(\alpha)(x\ot m)=\alpha(x)(1\ot m)$,  $\psi(\beta)(x)(h\ot m)=\beta(h\cdot x\ot m)$, for any $h\in H$, $x\in X$, $m\in M$, are well-defined maps, one the inverse of each other. It follows straightforward that, when $M=N$, this isomorphism is an algebra map. 
\epf

\subsection{Module functors} Given two $H$-comodule algebras $K, S$, we shall explicitly describe  the category of $\Rep(H)$-module functors between the associated module categories.

Under these hypothesis, we shall denote by ${}_S^H\Mo_K$ the category of finite-dimensional $(S,K)$-bimodules that are also left $H$-comodules, with comodule structure a morphism of $(S,K)$-bimodules.

\begin{prop}\label{prop1} Asumme $K, S$ are right $H$-simple left $H$-comodule algebras, and ${}_K\Mo, {}_S\Mo$ the corresponding $\Rep(H)$-module categories. There are equivalences
$$\Rex({}_K\Mo, {}_S\Mo)\simeq {}_S\Mo_K, \quad \Fun_{\Rep(H)}({}_K\Mo, {}_S\Mo)\simeq {}_S^H\Mo_K.$$
\end{prop}
\pf  We shall only explain the definition of the equivalences. For the complete proof see \cite[Prop.1.23]{AM}. The first equivalence is a consequence of a Theorem of Watts, see \cite{Wa}. The functor $\Phi: {}_S\Mo_K\to \Rex({}_K\Mo, {}_S\Mo)$, $\Phi(B)(M)=B\ot_K M$, is an equivalence of categories.

If $P\in {}_S\Mo_K$, define $F_P :{}_K\Mo\to {}_S\Mo$  the functor given by $F_P(M)=P\ot_K
M$. The correspondence $P\mapsto F_P$ is an equivalence of categories.

 If  $P\in {}_S^H\Mo_K$, $X\in\Rep(H)$, $M\in {}_K\Mo$ the functor $F_P$ has a  module structure  as follows 
$$c_{X,M}:P\ot_K
(X\ot_{\ku} M)\to X\ot_{\ku} (P\ot_K M),$$ 
$$c_{X,M}(p\ot  x\,\ot
m)=p\_{-1}\cdot x \,\ot\, p\_0\ot m, $$ 
for any $ p\in P, x\in X, m\in M.$ Here the map $\lambda:P\to H\ot P$, $\lambda(p)=p\_{-1} \ot p\_0$, is the left $H$-coaction of $P$.\epf
The next result is a direct consequence of Proposition \ref{prop1}.
\begin{cor}\label{dual:caseHopf} Let $K$ be a right $H$-simple left $H$-comodule algebra. There is a monoidal equivalence
$$\Rep(H)^*_{{}_K\Mo}\simeq {}_K^H\Mo_K.$$\qed
\end{cor}

Assume $K, S$ are $H$-comodule algebras. The category ${}_S\Mo_K$ has a $\Rep(H)$-bimodule structure as follows. If $P\in {}_S\Mo_K$ , $X\in \Rep(H)$, then 
$$X\otb P= X\otk P, \quad P\otb X= P\ot_K (X\otk K).$$
The $(S,K)$-action on the spaces $X\otb P, P\otb X$ are
$$s\cdot (x\ot p)\cdot k = s\_{-1}\cdot x\ot s\_0\cdot p\cdot k, $$
$$ s\cdot (p\ot (x\ot l))\cdot k =s\cdot p\ot  (x\ot lk),$$
for any $s\in S, k,l\in K, p\in P, x\in X$. The natural isomorphisms relating both actions are given by
\begin{equation}\label{bimod-iso-hopf} \gamma_{X,P,Y}: (X\otb P)\otb Y\to  X\otb ( P\otb Y), \gamma_{X,P,Y}((x\ot p)\ot (y\ot k))= x\ot (p\ot y\ot k),
\end{equation}
for any $X, Y\in \Rep(H)$, $P\in {}_S\Mo_K$, $x\in X, y\in Y, p\in P$, $k\in K$. It follows by a straightforward computation that the maps $\gamma_{X,P,Y}$ satisfy \eqref{bimod11}, \eqref{bimod12} and \eqref{bimod13}.

Recall that $\Rex({}_K\Mo, {}_S\Mo)$ has a $\Rep(H)$-bimodule category structure, see \eqref{bimod-action-funct}. The proof of the next lemma follows straightforward.
\begin{lema}\label{bimod-comodalg} The equivalence
$\Rex({}_K\Mo, {}_S\Mo)\simeq {}_S\Mo_K$ presented in  Proposition \ref{prop1} is an equivalence of $\Rep(H)$-bimodule categories.\qed
\end{lema}

\subsection{The center of dual tensor categories}\label{Section:center-dual-Hopf}

In Section \ref{Section:moritainvariance},  for any exact $\ca$-module category $\Mo$ we presented an equivalence of braided tensor categories $\theta_\Mo:\Zc(\ca)\to \Zc(\ca^*_\Mo)$. In this section, we shall explicitly give this equivalence in the case $\ca=\Rep(H)$ and $\Mo={}_K\Mo$ for a right $H$-simple left $H$-comodule algebra $K$. For this, we shall use the monoidal equivalences $\Zc(\Rep(H))\simeq {}^H_H\YD$, and  $\Rep(H)^*_{{}_K\Mo}\simeq {}_K^H\Mo_K$. The last one presented in Corollary \ref{dual:caseHopf}. 

Set $\theta_K=\theta_{{}_K\Mo}:{}^H_H\YD\to \Zc({}_K^H\Mo_K)$. If $V\in {}^H_H\YD$ then $\theta_K(V)= V\otk K$. The $K$-bimodule  and left $H$-comodule structure are given by
$$ k\cdot (v\ot t)\cdot s= k\_{-1}\cdot v\ot k\_0ts,$$
$$\lambda(v\ot t)= v \_{-1}t\_{-1}\ot v\_0\ot t\_0,$$
for any $v\in V$, $t,k,s\in K$. The half braiding of the object $ V\otk K$ is given by
$$\sigma^{V}_P:(V\otk K)\ot_K P\to P\ot_K (V\otk K),$$
$$\sigma^{V}_P(v\ot t\ot p)=(t\cdot p)\_0 \ot \Sc^{-1}((t\cdot p)\_{-1})\cdot v \ot 1$$
for any $P\in {}_K^H\Mo_K$, $v\in V, p\in P$, $t\in K$. This formula comes from \eqref{braiding-center}.
\section{ The character algebra for representations of $\Rep(H)$}\label{Section:character-algebramain}

Given a finite dimensional Hopf algebra $H$, and $\Mo$ a representation of the tensor category  $\Rep(H)$. We aim to compute the adjoint algebra $\cha_{\Mo}$ and the corresponding space of class functions as introduced by K. Shimizu \cite{Sh1}, \cite{Sh2}.

\subsection{The adjoint algebra  and the space of class functions }\label{SubSection:character algebra}

Let $\ca$ be a finite tensor category, and let $\Mo$ be an exact indecomposable left module category over $\ca$. We shall further assume that $\Mo$ is strict. First, we shall recall the definition of the algebra $\cha_\Mo\in \Zc(\ca)$.

The action functor $\rho_\Mo:\ca\to \Rex(\Mo)$ is
$$\rho_\Mo(X)(M)=X\otb M, \quad X\in \ca, M\in \Mo.$$
It was proven in \cite[Thm. 3.4]{Sh2} that the right adjoint of $\rho_\Mo$ is the functor $\rho_\Mo^{\ra}:  \Rex(\Mo)\to\ca$, such that for any $F\in \Rex(\Mo)$
$$ \rho_\Mo^{\ra}(F)= \int_{M\in \Mo} \uhom(M, F(M)).$$

The counit and unit of the adjunction $(\rho_\Mo, \rho_\Mo^{\ra})$, will be denoted by 
$$\epsilon:\rho_\Mo\circ  \rho_\Mo^{\ra}\to \Id_{\Rex(\Mo)},\quad \eta:\Id_\ca\to \rho_\Mo^{\ra}\circ \rho_\Mo.
$$

According to Lemma \ref{modfunct-adjoint} the functor $\rho_\Mo^{\ra}$ has a structure of $\ca$-bimodule functor  as follows. 
The left and right module structure of $\rho_\Mo^{\ra}$ are
\begin{align}\label{xil}\begin{split}
\xi^l_{X,F}: \rho_\Mo^{\ra}(X\otb F)\to X\otb\rho_\Mo^{\ra}(F),\\
\big(\xi^l_{X,F}\big)^{-1} = \rho_\Mo^{\ra}(\id_X\otb \epsilon_F) \eta_{X\otb \rho_\Mo^{\ra}(F)},
\end{split}
\end{align}
\begin{align}\label{xir}\begin{split}
\xi^r_{X,F}: \rho_\Mo^{\ra}(F\otb X)\to
\rho_\Mo^{\ra}(F)\otb X,\\
\big(\xi^r_{X,F}\big)^{-1}=\rho_\Mo^{\ra}(\epsilon_F\otb \id_X)  \eta_{\rho_\Mo^{\ra}(F)\otb X},
\end{split}
\end{align}
for any $X\in \ca$, $F\in \Rex(\Mo)$.  This description appears in \cite[Equation A.9]{Sh2}.

Since the functor $\rho_\Mo^{\ra}:  \Rex(\Mo)\to\ca$ is a $\ca$-bimodule functor, we can consider the  functor $\Zc_\ca(\rho_\Mo^{\ra}): \End_\ca(\Mo)\to \Zc(\ca)$. Here $\Zc_\ca$ is the 2-functor described in Section \ref{Section:relative center}.

\begin{defi}\cite[Subsection 4.2]{Sh2} The \textit{adjoint algebra} of the module category $\Mo$ is the algebra in the center of $\ca$, $\cha_\Mo:=\Zc(\rho_\Mo^{\ra})(\Id_\Mo)\in \Zc(\ca)$. The \textit{adjoint algebra of the tensor category }$\ca$ is the algebra $\cha_\ca$ of the regular module category $\ca$.
\end{defi}

It was explained in \cite[Subection 4.2]{Sh2} that the algebra structure of $\cha_\Mo$ is given as follows. Let $\pi_\Mo: \cha_\Mo\xrightarrow{ .. } \uhom( -,-)$ denote the dinatural transformation of the end $\cha_\Mo$. The product and the unit of $\cha_\Mo$ are 
$$m_\Mo:\cha_\Mo\ot \cha_\Mo\to \cha_\Mo, \quad u_\Mo:\uno\to \cha_\Mo,$$
defined to be the unique morphisms such that they satisfy
\begin{equation}\label{product-cha-M}
\begin{split}\pi_\Mo(M)\circ m_\Mo= \comp^\Mo_{M}\circ (\pi_\Mo(M)\ot \pi_\Mo(M)),\\ \pi_\Mo(M)\circ  u_\Mo= \coev^\Mo_{\uno, M},\end{split}
\end{equation}

for any $M\in \Mo$. For the definition of $\coev^\Mo$ and $ \comp^\Mo$ see Section \ref{subsection:internal hom}.

\begin{defi}\label{definition:classfunct}\cite[Definition 5.1]{Sh2} The \textit{space of class functions} of $\Mo$ is $\cf(\Mo):=\Hom_\ca(\Fc_\ca(\cha_\Mo),\uno)=\Hom_{\Zc(\ca)}( \cha_\Mo,  \cha_\ca)$.
\end{defi}
The following  result will be useful when computing the adjoint algebra in particular examples. The first three statements are contained in \cite{Sh1}, \cite{Sh2}.
\begin{lema}\label{fpdim-char} Let $\Mo$ be an exact indecomposable $\ca$-module category. The following statements hold.
\begin{itemize}
\item[(i)] If $I:\ca\to \Zc(\ca) $ is a right adjoint to the forgetful functor $\Fc:\Zc(\ca)\to \ca$, then $\cha_{\ca} \simeq I (\uno)$.
\item[(ii)] There exists an isomorphism $\cha_{\ca^*_\Mo}\simeq \theta_\Mo(\cha_\Mo)$ as algebra objects in  $\Zc(\ca^*_\Mo)$.
\item[(iii)] $\cf(\ca^*_\Mo)\simeq \End_{\Zc(\ca)}(\cha_\Mo)$.
\item[(iv)] $\fp(\cha_\Mo)=\fp(\ca)$.
\end{itemize}
\end{lema}
\pf Recall that $\theta_\Mo:\Zc(\ca)\to \Zc(\ca^*_\Mo)$ is the braided equivalence presented in Section \ref{Section:moritainvariance}.

It is proven in \cite[Corollary 3.15]{Sh2}  that the functor $\theta_\Mo \circ \Zc(\rho_\Mo^{\ra})$ is the right adjoint of the forgetful functor $\Fc_{\ca^*_\Mo}$. Taking $\Mo=\ca$, this implies part (i), and taking arbitrary $\Mo$ follows part (ii).

\medbreak
(iii). This is \cite[Theorem 5.12]{Sh2}.
\medbreak
(iv). Let $F:\Zc(\ca)\to \ca$ be the forgetful functor, and $I:\ca\to \Zc(\ca)$ its right adjoint. It was proven in \cite[Proposition 7.16.5 ]{EGNO} that $\fp(I(\uno))=\fp(\ca)$. Hence $\fp(\cha_\ca)=\fp(\ca)$, for any finite tensor category $\ca$. Applying this result to $\ca^*_\Mo$ we obtain that 
$$\fp(\cha_\Mo)=\fp(\cha_{\ca^*_\Mo})=\fp(\ca^*_\Mo)=\fp(\ca).$$
The first equality follows from part (ii), and the last  one is \cite[Corollary 3.43]{EO}.
\epf

\subsection{The adjoint algebra for module categories over Hopf algebras}

Let $H$ be a finite dimensional Hopf algebra. Let $K$ be a finite-dimensional left $H$-comodule algebra. The category ${}_K\Mo$ is a left $\Rep(H)$-module category. See Section \ref{subsection:modcat-over-hopf}. We aim to compute $\cha_K=\cha_{{}_K\Mo}$ as an algebra in the category ${}^H_H\YD$ of Yetter-Drinfeld modules over $H$.  For this, we shall explicitly give a description of  the functor $\Zc(\rho_\Mo^{\ra})$.

Identifying $\Rex({}_K\Mo)={}_K\Mo_K$, we shall denote by $\rho_K:\Rep(H)\to  {}_K\Mo_K$, the action functor. Explicitly, if $X\in \Rep(H)$ then 
$$\rho_K(X)=X\otk K.$$ The left and right $K$-action on 
$X\otk K$ are given by:
$$ s\cdot (x\ot k)\cdot t=s\_{-1}\cdot x \ot s\_0kt, $$
for any $x\in X, s,t,k\in K$. 
\begin{defi} For any $P\in  {}_K^H\Mo_K$, define  $S^K(H,P)$ as the space of left $K$-linear morphisms $\alpha\in \Hom_K(H\otk K,P)$ such that for any $k\in K, h\in H$
\begin{equation}\label{coend-hopfc}
    \alpha(h\ot k )=\alpha(h\ot  1)\cdot k.
\end{equation}
\end{defi}
The space $S^K(H,P)$ has a left $H$-module structure $\cdot:H\otk S^K(H,P)\to S^K(H,P)$ and a left $H$-comodule structure $\lambda:S^K(H,P)\to H\otk S^K(H,P)$,  defined  by:
\begin{equation}\label{H-mod-coend}
(h\cdot \alpha)(x\ot k)=\alpha(xh\ot k),
\end{equation}
\begin{equation}\label{H-comod-coend}
\lambda: S^K(H,P)\to H\otk S^K(H,P),\quad \lambda(\alpha)=\alpha^{-1}\ot \alpha^0,
\end{equation}
for any $\alpha\in S^K(H,P)$, $h, x\in H, k\in K$. Here for any $h\in H, k\in K$
\begin{equation}\label{coact-coend-def}
\alpha^{-1}\ot \alpha^0(h\ot k)=\Sc(h\_1) \alpha(h\_2\ot 1)\_{-1} h\_3\ot \alpha(h\_2\ot 1)\_0 \cdot k.
\end{equation}
When $P=K$, we shall denote $S(H,K):=S^K(H,K)$. It follows by a  straightforward computation that,    \eqref{H-mod-coend}, \eqref{H-comod-coend} are well defined maps, and they define an $H$-action and a $H$-coaction.

\begin{lema} The space $S^K(H,P)$ is an object in the category ${}^H_H\YD$.

\end{lema}
\pf 
We must prove compatibility condition \eqref{yd-module},  that is
\begin{equation}\label{yd:inS}
\lambda(x\cdot \alpha)=x_{(1)} \alpha^{-1} S(x_{(3)}) \ot x_{(2)} \cdot \alpha^0,
\end{equation}
for any $x \in H$ and for all $\alpha \in S^K(H,P)$. Take $\phi \in H^*$, $h, x \in H$ and $k \in K$. Evaluating the right hand side of \eqref{yd:inS} in $\phi\ot h\ot k$ gives
\begin{align*}
    \big< \phi, & x_{(1)} \alpha^{-1} S(x_{(3)}) \big>\, (x_{(2)} \cdot \alpha^0)(h \ot k)= \\ 
    &= \big<\phi\_1, x_{(1)} \big> \big<\phi\_3, S(x_{(3)}) \big> \big<\phi\_2, \alpha^{-1} \big> x_{(2)} \cdot \alpha^0(h \ot k)\\ 
    &= \big< \phi\_1, x_{(1)} \big> \big< \phi\_3, S(x_{(3)}) \big> \big< \phi\_2, \alpha^{-1} \big> \alpha^0(hx_{(2)}\ot k) \\ 
    &= \big< \phi\_1, x_{(1)} \big> \big< \phi\_3, S(x_{(3)}) \big> \\ & \,\,\big<\phi\_2, S((hx_{(2)})_{(1)}) \alpha((hx_{(2)})_{(2)} \ot 1)_{(-1)} (hx_{(2)})_{(3)} \big> 
    \alpha((hx_{(2)})_{(2)} \ot 1)\_0 \cdot k \\
    &= \big< \phi, x\_1S(x\_2) S(h\_1)\alpha(h\_2x\_3\ot 1)_{(-1)}h\_3x\_4S(x\_5)\big>\\
  &\quad   \quad \alpha(h\_2x\_3\ot 1)\_0 \cdot k \\
    &= \big< \phi, S(h_{(1)}) \alpha(h_{(2)}x\ot 1)_{(-1)} h_{(3)} \big> \,\alpha(h_{(2)}x \ot 1)\_0 \cdot k \\
    &= \big< \phi, S(h_{(1)}) (x \cdot \alpha)(h_{(2)} \ot 1)_{(-1)} h_{(3)}\big> (x \cdot \alpha)(h_{(2)} \ot 1)\_0 \cdot k \\
    &= \big<\phi, (x \cdot \alpha)^{-1} \big> (x \cdot \alpha)^0 (h \ot k).
\end{align*}\epf

\begin{lema} The space $S(H,K)$ is identified with the subspace of elements $\sum_i f_i\ot k_i\in H^*\otk K$ such that for any $t\in K$, $h\in H$
\begin{equation}\label{sk-elements} \sum_i <f_i, t\_{-1}h> k_i t\_0= \sum_i <f_i, h> tk_i.
\end{equation}
\end{lema}
\pf We explained in Section \ref{Section:homint-hf} that the space $\Hom_K(H\otk K, K)$ can be identified with elements  $\sum_i f_i\ot T_i\in H^*\otk \End(K)$ such that they verify \eqref{internal-id}. An element $\sum_i f_i\ot T_i\in H^*\otk \End(K)$, with  $\{f_i\}_i$  linearly independent, belongs to $S(H,K)$ if it satisfies \eqref{coend-hopfc}. This means that $T_i(k)=T_i(1)k$, for all $i$, and all $k\in K$. Thus, we can identify $T_i$ with the left multiplication map by the element   $T_i(1)$. Under this identification \eqref{internal-id} is equivalent to \eqref{sk-elements}.
\epf

\begin{rmk} If $K\subseteq H$ is a left coideal subalgebra, then  elements of the space $(H/K^+H)^*\otk Z(K)$ are inside of $S(H,K)$. Here $Z(K)$ is the center of $K$. This observation follows from the fact that $(H/K^+H)^*$ can be identified with elements $f\in H^*$ such that $<f,kh>=<\epsilon, k> <f,h>$, for any $k\in K$, $h\in H$.
\end{rmk}

\begin{teo}\label{char-for-mod-hopf} Let $K$ be a finite dimensional left $H$-comodule algebra and $P\in  {}_K^H\Mo_K$. Then $\Mo={}_K\Mo$ is a left $\Rep(H)$-module category. There is an isomorphism of $H$-modules
$$S^K(H,P)\simeq \int_{M\in \Mo} \uhom(M,P \ot_K M).$$
When $P=K$, this isomorphism is an algebra map.
\end{teo}
\pf We shall use the description of the internal Hom of the module category ${}_K\Mo$ given in Lemma \ref{int-hom-hopf}.  First, we shall prove that $S^K(H,P)= \int_{M\in \Mo} \uhom(M,P \ot_K M)$ as objects in $\Rep(H)$. Observe that if $M, M', N, N'$ are objects in ${}_K\Mo$, and $f:M\to M'$, $g:N\to N'$ are $K$-module morphisms, then the functor $\uhom:\Mo^{\op}\times \Mo\to \Rep(H)$ is defined on morphisms by
$$\uhom(f,g):\Hom_K(H\otk M', N)\to \Hom_K(H\otk M, N'),$$
$$\alpha\mapsto g\circ \alpha\circ (\id_H\ot f).$$
For any $M\in {}_K\Mo$ define $\pi^P_M:S^K(H,P)\to \Hom_K(H\otk M, P \ot_K M)$, by 
$$ \pi^P_M(\alpha)(h\ot m)=\alpha(h\ot 1)\ot_K m,$$ for any $h\in H$, $m\in M$. Equation \eqref{coend-hopfc} implies that $\pi^P_M(\alpha)$ is a $K$-module morphism. It follows straightforward that $\pi^P_M$ is an $H$-module map and that it is dinatural.  

Assume that $(E,d)$ is a pair, where $E\in \Rep(H)$, and $$d:E\xrightarrow{..} \uhom(-,  P\ot_K -)$$ is a dinatural transformation. Dinaturality, in this case, implies that
for any pair of $K$-modules $M, N$, and a $K$-module map $f:M\to N$, we have
\begin{equation}\label{dinat-hopf-case}
(\id_P\ot f)\circ d_M(e)=d_N(e) (\id_H\ot f),
\end{equation}
for any $e\in E$. In particular, if $N$ is any $K$-module, and $n\in N$, define $f_n:K\to N$, $f_n(k)=k\cdot n$. Hence $f_n$ is a $K$-module map, whence, equation \eqref{dinat-hopf-case} implies that $(\id_P\ot f_n)\circ d_K(e)=d_N(e) (\id_H\ot f_n)$.  Evaluating this equality  in the element $h\ot 1\in H\otk K$ we have that
\begin{equation}\label{dinat-hopf-case2} d_K(e)(h\ot 1)\ot n=d_N(e)(h\ot n), \end{equation}
for any $h\in H$. This implies, taking $N=K$, that the element $d_K(e)\in S^K(H,P)$, that is, $d_K(e)$ satisfies \eqref{coend-hopfc}.

Define $\phi:E\to S^K(H,P)$ as $\phi=d_K$. Then, equation \eqref{dinat-hopf-case2} implies that $d_N= \pi^P_N\circ \phi$ for any $K$-module $N$.  This proves that the object $S^K(H,P)$ together with the dinatural transformations $\pi^P$ satisfies the universal property of the end. Thus $S^K(H,P)= \int_{M\in \Mo} \uhom(M,P\ot_K M).$ 

When $P=K$, it is not difficult to verify that the product of the adjoint algebra defined in terms of the dinatural transformation, see  \eqref{product-cha-M}, coincides with the product described in \eqref{product-internal-hom}.
\epf

So far, we have described the structure of the end  $\int_{M\in \Mo} \uhom(M,P \ot_K M)$ as an object in $\Rep(H)$. It remains to describe the structure as an object in the category of Yetter-Drinfeld modules over $H$. 
The next two results will be initial steps towards this objective.
\medbreak
Define the functor $\rhob_K: {}_K\Mo_K\to \Rep(H),$ $\rhob_K(P)=S^K(H,P)$ for any $P\in  {}_K\Mo_K.$ If $P, Q\in {}_K\Mo_K$ and $f:P\to Q$ is a morphism of $K$-bimodules, then
$$\rhob_K(f): S^K(H,P)\to S^K(H,Q),\quad \rhob_K(f)(\alpha)=f\circ \alpha.$$
\begin{prop}\label{adjoint-rhoK} The functor $\rhob_K: {}_K\Mo_K\to \Rep(H)$ is the right adjoint of the functor $\rho_K$. The unit and counit of the adjunction $(\rho_K, \rhob_K)$ are given by
$$\eta:\Id_{\Rep(H)}\to \rhob_K\circ \rho_K,\quad \epsilon:\rho_K\circ  \rhob_K\to \Id_{{}_K\Mo_K} ,$$
$$\eta_X(x)(h\ot k)=h\cdot x\ot k, \quad\epsilon_P(\alpha\ot k)=\alpha(1\ot k),$$
for any $X\in \Rep(H), P\in {}_K\Mo_K$, $x\in X, h\in H$, $\alpha\in S^K(H,P), k\in K$.
\end{prop}

\pf For any $X\in \Rep(H), P\in {}_K\Mo_K$, $x\in X, k\in K$, $h\in H$ define
\begin{equation}\label{adjoint-rho1} \begin{split}\phi_{X,P}:&\Hom_H(X, S^K(H,P))\to \Hom_{(K,K)}(X\otk K,  P),\\
&\phi_{X,P}(\alpha)(x\ot k)=\alpha(x)(1\ot k),
\end{split}
\end{equation}
\begin{equation}\label{adjoint-rho2}\begin{split} \psi_{X,P}:&\Hom_{(K,K)}(X\otk K,  P)\to\Hom_H(X, S^K(H,P)),\\
&\psi_{X,P}(\beta)(x)(h\ot k)=\beta(h\cdot x\ot k).
\end{split}
\end{equation}
It follows by a straightforward computation that the maps $\phi_{X,P}, \psi_{X,P}$ are natural morphisms, 
and are inverses of each other.  The unit and counit of the adjunction are given by
$$\eta_X= \psi_{X,X\otk K}(\id_{X\otk K}), \quad \epsilon_P=\phi_{S^K(H,P),P}(\id_{S^K(H,P)}).$$\epf

The next result is a particular case of Example \ref{relat-center-rex}. Since we need the explicit equivalence, we shall write the proof.
\begin{lema}\label{equiv-relat-c} There is an equivalence of categories 
$\Zc_{\Rep(H)}({}_K\Mo_K)\simeq  {}_K^H\Mo_K$.
\end{lema}
\pf Let $(M, \sigma)\in \Zc_{\Rep(H)}({}_K\Mo_K)$. This means that $M\in {}_K\Mo_K$, and  the half-braiding is given by $\sigma_X: M\ot_K (X\otk K)\to X\otk M$, for any $X\in \Rep(H)$. Define $\lambda:M\to H\otk M$, $\lambda(m)= \sigma_H(m\ot 1_H\ot 1_K)$ for any $m\in M$. This establishes a functor 
$$\Phi: \Zc_{\Rep(H)}({}_K\Mo_K)\to {}_K^H\Mo_K, \quad \Phi(M, \sigma)=(M,\lambda).$$

If $(M,\lambda)\in {}_K^H\Mo_K$, define $\sigma^\lambda_X: M\ot_K (X\otk K)\to X\otk M$, the map
$$\sigma^\lambda_X(m\ot x\ot k)=m\_{-1}\cdot x\ot m\_0\cdot k, $$
for any $X\in \Rep(H)$, $m\in M$, $k\in K$. It follows by a simple computation that $\sigma^\lambda_X$ is a well-defined isomorphism, it is a $K$-bimodule map and it satisfies \eqref{half-braid}. This defines a functor $\Psi: {}_K^H\Mo_K\to \Zc_{\Rep(H)}({}_K\Mo_K)$, $\Psi(M,\lambda)=(M, \sigma^\lambda)$. 
\epf

For any  $P\in {}_K^H\Mo_K$ recall the structure of Yetter-Drinfeld module over $H$ of $S^K(H,P)$ given by \eqref{H-mod-coend}, \eqref{H-comod-coend}.
\begin{teo}\label{ydd-ra} For any $P\in {}_K^H\Mo_K$ there is an isomorphism $S^K(H,P)\simeq \Zc(\rhob_K)(P)$ as objects in ${}^H_H\YD$.
\end{teo}
\pf
If $P\in {}_K^H\Mo_K$, then $\Zc(\rhob_K)(P,\sigma^\lambda)=(S^K(H,P), \sigma^P),$ where, according to \eqref{relative-bimod-funct}, the half braiding for $\rhob_K(P)=S^K(H,P)$ is the morphism $\sigma^P_X: S^K(H,P)\otk X \to X\otk S^K(H,P)$, given by the composition
$$S^K(H,P)\otk X\xrightarrow{ (\xi^r_{X, P})^{-1}} S^K(H,P\ot_K (X\otk K))\xrightarrow{ \rhob_K(\sigma^\lambda_X)} $$
$$\to S^K(H,X\otk P) \xrightarrow{ \xi^l_{X, P}} X\otk  S^K(H,P), $$
for any $X\in \Rep(H)$. Recall that $\sigma^\lambda$ is the half braiding associated to $P$ explained in Lemma \ref{equiv-relat-c}. To compute $\sigma^P$, we need to compute the bimodule structure of the functor $\rhob_K$. Both structures are given by equations \eqref{xil}, \eqref{xir}.  

Using  the formula for the unit and counit of the adjunction $(\rho_K, \rhob_K)$ given in Proposition \ref{adjoint-rhoK} we obtain that 
\begin{align*} ( \xi^l_{X, P})^{-1}(x\ot \alpha)(h\ot k)&=\rhob_K(\id_X\ot \epsilon_P)\eta_{X\otk S^K(H, P)}(x\ot \alpha)(h\ot k)\\
&=(\id_X\ot \epsilon_P)\eta_{X\otk S^K(H, P)}(x\ot \alpha)(h\ot k)\\
&=(\id_X\ot \epsilon_P)(h\_1\cdot x\ot h\_2\cdot \alpha\ot k)\\
&=h\_1\cdot x\ot   h\_2\cdot \alpha(1\ot k)\\
&=h\_1\cdot x\ot \alpha(h\_2\ot k),
\end{align*}
and
\begin{align*}( \xi^r_{X, P})^{-1}(\alpha\ot x)(h\ot k)&= (\epsilon_P\ot\id_{X\otk K})\eta_{S^K(H, P) \otk X} (\alpha\ot x)(h\ot k)\\
&=(\epsilon_P\ot\id_{X\otk K}) (h\_1\cdot \alpha\ot 1\ot h\_2\cdot x\ot k) \\
&=(h\_1\cdot \alpha)(1\ot 1)\ot  h\_2\cdot x\ot k\\
&=\alpha(h\_1\ot 1)\ot  h\_2\cdot x\ot k,
\end{align*}
for any $\alpha\in S^K(H, P)$, $x\in X$, $h\in H$, $k\in K$.

Now, the $H$-coaction of $\Zc(\rhob_K)(P,\sigma^\lambda)$ associated with the half braiding $\sigma^P$ is 
$$\lambda^P:S^K(H, P)\to H\otk S^K(H, P), \quad \lambda^P(\alpha)=\sigma^P_H(\alpha\ot 1_H).$$
Let us denote $\lambda^P(\alpha)=\alpha^{-1}\ot \alpha^0$. Using the formula for $\sigma^P_X$, we know that $( \xi^l_{H, P})^{-1}\sigma^P_H(\alpha\ot 1_H)=\rhob_K(\sigma^\lambda_H) ( \xi^r_{H, P})^{-1}(\alpha\ot 1_H)$. Evaluating this equality in $h\ot k\in H\otk K$ we obtain that
\begin{align*} ( \xi^l_{H, P})^{-1}\sigma^P_H(\alpha\ot 1_H)(h\ot k)&=h\_1\alpha^{-1}\ot \alpha^0(h\_2\ot k),
\end{align*}
is equal to
\begin{align*} \rhob_K(\sigma^\lambda_H)  \xi^r_{H, P}(\alpha\ot 1_H)(h\ot k)&= \sigma^\lambda_H( \xi^r_{H, P})^{-1}(\alpha\ot 1_H)(h\ot k)\\
&=\sigma^\lambda_H(\alpha(h\_1\ot 1) \ot h\_2\ot k)\\
&=\alpha(h\_1\ot 1)\_{-1}h\_2\ot \alpha(h\_1\ot 1)\_0\cdot k.
\end{align*}
Thus 
\begin{equation}\label{th11} h\_1 \alpha^{-1}\ot \alpha^0(h\_2\ot k)=\alpha(h\_1\ot 1)\_{-1}h\_2\ot \alpha(h\_1\ot 1)\_0\cdot k,
\end{equation}
for any $\alpha\in S^K(H, P)$, $h\ot k\in H\otk K$. Hence
\begin{align*}\alpha^{-1}\ot \alpha^0(h\ot k)&=\alpha^{-1}\ot \alpha^0(\epsilon_H(h\_1)h\_2\ot k)\\
&= \Sc(h\_1)h\_2\alpha^{-1}\ot \alpha^0(h\_3\ot k)\\
&=\Sc(h\_1) \alpha(h\_2\ot 1)\_{-1} h\_3\ot \alpha(h\_2\ot 1)\_0 \cdot k.
\end{align*}
The last equality follows from \eqref{th11}. This formula coincides with \eqref{coact-coend-def}.
\epf
As a consequence of Theorems \ref{char-for-mod-hopf} and \ref{ydd-ra} we obtain the next result.
\begin{cor} Let $K$ be a finite dimensional right $H$-simple left $H$-comodule with trivial coinvariants. There exists an isomorphism of algebras 
$$ S(H,K)\simeq \cha_K,$$
in the category ${}^H_H\YD$.\qed
\end{cor}

\begin{exa}[Case $K=H$] We denote by $H_{ad}$ the algebra in the category ${}^H_H\YD$ whose underlying algebra is $H$, with $H$-coaction given by the coproduct and  $H$-action given by the adjoint action, that is $h\triangleright x=h\_1 xS(h\_2)$, $h,x\in H$. Since $H$ is an $H$-comodule algebra with the coproduct, we can consider $S^H(H,H)$. The map $\phi:S(H,H)\to H_{ad}$, $\phi(\alpha)=\alpha(1\ot 1)$ is an isomorphism of algebras in  ${}^H_H\YD$. Indeed, it is an $H$-module map. Take $\alpha\in S(H,H)$, $h, t\in H$, then
\begin{align}\label{iso1}\begin{split} \alpha(h\ot t)&=\alpha(h_1\ot h_2\Sc(h_3)t)=h\_1 \alpha(1\ot \Sc(h_2)t)\\
&= h\_1 \alpha(1\ot 1)\Sc(h_2)t.
\end{split}
\end{align}
The second equality because $\alpha$ is an $H$-module map, and the last equality follows from \eqref{coend-hopfc}. Then
\begin{align*}\phi(h\cdot \alpha)&=(h\cdot \alpha)(1\ot 1)=
\alpha(h\ot 1)\\
&= h\_1 \alpha(1\ot 1)\Sc(h_2)= h\triangleright \phi(\alpha).
\end{align*}
\end{exa}
It follows by a straightforward computation that, $\phi$ is an  algebra and an $H$-comodule map. Using \eqref{iso1}, it follows that the map $\psi: H_{ad}\to S(H,H)$, $\psi(x)(h\ot t)=h\_1 x\Sc(h_2)t$, $x,h,t\in H$, is the inverse of $\phi$.
\begin{exa}[Case $K=\ku$] We denote by $H^*_{ad}$ the following algebra in the category ${}^H_H\YD$.   The underlying algebra is $H^*$. The $H$-action and $H$-coaction are $\cdot : H\otk H^*_{ad}\to H^*_{ad}$, $\lambda: H^*_{ad}\to H\otk H^*_{ad}$, $\lambda(f)=f\_{-1}\ot f\_0$, where 
$$(h\cdot f)(x)=f(xh), \quad <g, f\_{-1}>  f\_0= \Sc(g\_1)fg\_2,$$
for any $h,x\in H$, $g\in H^*$. It follows that $S(H,\ku)=H^*_{ad}$.
\end{exa}

\section{Some explicit calculations}\label{Section:example-comp}

In this section we shall explicitly compute the adjoint algebra for the representation categories  of group algebras and their duals.
We shall use  the identification of $S(H,K)$ with elements in $H^*\otk K$ such that they satisfy \eqref{sk-elements}. First we recall the classification of exact indecomposable module categories over group algebras and their duals.

\subsection{Module categories over the tensor categories $\Rep(\ku^G)$, $\Rep(\ku G)$}\label{Section:modgroup} Assume $G$ is a finite group. We shall recall the classification of exact indecomposable module categories over $\Rep(\ku^G)$ and $\Rep(\ku G)$. For this, we shall give families of simple left $H$-comodule algebras, where $H=\ku^G, \ku G$.

Assume  $F\subseteq G$ is a subgroup and  $\psi\in Z^2(F,\ku^{\times})$ a 2-cocycle. We denote by $\ku_\psi F$ the twisted group algebra.  We can choose $\psi$ (in a cohomology class) such that 
$$\psi(f,g)\psi(g^{-1}, f^{-1})=1, \quad \psi(f,1)=\psi(1,f)=1, $$ for any $f,g\in F$. In such case we shall say that $\psi$ is \textit{normalized}.

The twisted group algebra $\ku_\psi F$ is a left $\ku G$-comodule algebra as follows. Elements in $\ku_\psi F$ are linear combinations of $e_f, f\in F$. The product and left $\ku G$-coaction are
$$ e_f e_h=\psi(f,h)\; e_{fh},\quad \lambda(e_f)=f\ot e_f,$$
for any $f,h\in F$. 
If $V$ is a simple $\ku_\psi F$-module, we can form the following algebra. The endomorphism algebra $\End(V)$ is a right $\ku F$-module, with action given by
$$(T\cdot  f)(v)=f^{-1}\cdot T(f\cdot v),$$
for any $f\in F, v\in V, T\in \End(V)$. Define $\kc(F,\psi, V)= \End(V)\ot_{\ku F} \ku G$. Let $S\subseteq G$ be a set of representative elements of cosets $F\backslash G$. 

Any element in $\kc(F,\psi, V)$ is of the form $\overline{T\ot s}$, for some $s\in S, T\in \End(V)$. Here $\overline{z}$ denotes the class of the element $z\in \ku G\otk \End(V)$ in the quotient $\ku G\ot_{\ku F} \End(V)$. The product in $\kc(F,\psi, V)$ is defined as follows:
$$(\overline{ T\ot x})(\overline{ U\ot y})=\delta_{x,y}\; \overline{T\circ U\ot x},$$
for any $T,U\in \End(V), x,y\in S$. The unit is $\sum_{s\in S} \overline{\Id\ot s}$.
The vector space $\kc(F,\psi, V)$ has a structure of right $\ku G$-module that makes it into a module algebra. The right action is:
$$(\overline{ T\ot x})\cdot g= \overline{ T\ot xg}, \quad g\in G.$$
With this action $\kc(F,\psi, V)$ is a right $\ku G$-module algebra, hence it is a left $\ku^G$-comodule algebra, with coaction
$$\lambda:\kc(F,\psi, V)\to \ku^G\otk \kc(F,\psi, V), \lambda(k)=k\_{-1}\ot k\_0,$$
such that for any $g\in G$, $<k\_{-1}, g> k\_0=k\cdot g$.

The next result is part of the folklore of representations of tensor categories. See, for example, \cite[Proposition 4.1, Lemma 4.3]{EO}. It also follow from Theorem \ref{sub-coid}.

\begin{teo}\label{mod-over-group} Let $G$ be a finite group. 
\begin{itemize}
\item[(i)] If $\Mo$ is an exact indecomposable module category over $\Rep(\ku G)$, there exists 
a subgroup $F\subseteq G$, a normalized 2-cocycle $\psi\in Z^2(F,\ku^{\times})$ such that $\Mo\simeq {}_{\ku_\psi F}\Mo$ as module categories.
\item[(ii)]  If $\Mo$ is an exact indecomposable module category over $\Rep(\ku^G)$, there exists 
a subgroup $F\subseteq G$, a normalized 2-cocycle $\psi\in Z^2(F,\ku^{\times})$, and a simple $\ku_\psi F$-module $V$ such that $\Mo\simeq {}_{\kc(F,\psi, V)}\Mo$ as module categories.\qed
\end{itemize}
\end{teo}

\begin{rmk} The equivalence class of the module category ${}_{\kc(F,\psi, V)}\Mo$ does not depend on the choice of the simple $\ku_\psi F$-module $V$. The twisted group algebra $\ku_\psi F$ is an algebra in the category $\Rep(\ku^G)$. One can prove that, regardless the choice of $V$, the module category ${}_{\kc(F,\psi, V)}\Mo$ is equivalent to $\Rep(\ku^G)_{\ku_\psi F}$.
\end{rmk}

\begin{rmk}\label{regular-modcat}  Let $F=\{1\}$ be the trivial subgroup of $G$, $\psi=1$ and $V=\ku$ with the trivial action. Denote $K=\kc(F,\psi, V)$. It is not difficult to see that $K\simeq \ku^G$ as left $\ku^G$-comodule algebras. Hence ${}_{\kc(\{1\},1, \ku)}\Mo\simeq \Rep(\ku^G)$ as $\Rep(\ku^G)$-module categories.
\end{rmk}

\subsection{Case $H=\ku G$}\label{Section:group-case} Let $F\subset G $ be a subgroup, and a normalized 2-cocycle $\psi\in Z^2(F,\ku^{\times})$. Let $K=\ku_\psi F$ be the twisted group algebra. We shall denote by $\{e_f\}_{f\in F}$ the canonical basis of $\ku_\psi F$. The product in this algebra is then $e_f e_l=\psi(f,l) e_{fl}$, for any $f, l\in F$.

Let $S\subseteq G$ be a set of representatives of right cosets $F\backslash G$ such that $1\in S$. 
Define $b:F\times F\to \ku^{\times}$ as
$$ b(l,f)=\frac{\psi(l,l^{-1}fl)}{\psi(f,l)}.$$
Also, for any $l\in F$, set $C_l=\{(g,f)\in F\times F: g^{-1}fg=l\}$.
For any $s\in S$, $l\in F$ define
$$\alpha_{s,l}=\sum_{(g,f)\in C_l} b(g,f)\,\delta_{gs}\ot e_f\in \ku^G\otk K.$$
Using the  identification explained in Subsection \ref{Section:homint-hf}, the element $\alpha_{s,l}$ can be seen as an element in $\Hom(\ku G\otk K, K)$, where 
\begin{equation}\label{morf-alphas} \alpha_{s,l}(x\ot e_h)=  \delta_{s,t} \,b(f,flf^{-1}) \psi(flf^{-1},h)\, e_{flf^{-1}h},
\end{equation}
if $x=ft\in G$, $t\in S$, $h\in F$.
\begin{lema}\label{basis} The set $B= \{\alpha_{s,l}\in\ku^G\otk \ku_\psi F: s\in S, l\in F\}$ is a basis of $S(\ku G, \ku_\psi F)$.
\end{lema}
\pf Clearly $B$ is a set of linearly independent elements. Let $z$ be an arbitrary element of $\ku^G\otk \ku_\psi F$. Thus $z=\sum_{x\in G, f\in F} \xi_{x,f}\, \delta_x\ot e_f$, for certain scalars $\xi_{x,f}\in \ku$. If $z\in S(\ku G, \ku_\psi F)$, equation \eqref{sk-elements} implies that
$$\sum_{x\in G, f\in F}  \xi_{x,f} \delta_x(ly) e_f e_l=\sum_{x\in G, f\in F}  \xi_{x,f} \delta_x(y) e_l e_f$$
 for any $y\in G, l\in F$. This implies that
 $$ \sum_{ f\in F}  \xi_{ly,f} \, \psi(f,l) \,e_{fl}=\sum_{ f\in F}  \xi_{y,f}\, \psi(l,f) \,e_{lf}. $$
This equality implies, by looking at the coefficient of $e_{lf}$, that
$$\xi_{ly,f}=\xi_{y,l^{-1}fl}\, b(l,f),$$
for any $l,f\in F$, $y\in G$. Whence  
\begin{align*} z&=\sum_{x\in G, f\in F} \xi_{x,f}\, \delta_x\ot e_f=\sum_{s\in S, g,f\in F} \xi_{gs,f}\, \delta_{gs}\ot e_f\\
&=\sum_{s\in S, g,f\in F} \xi_{s,g^{-1}fg}\, b(g,f)\, \delta_{gs}\ot e_f \\
&=\sum_{s\in S, l\in F} \sum_{(g,f)\in C_l}    \xi_{s,l}\, b(g,f)\, \delta_{gs}\ot e_f= \sum_{s\in S, l\in F}\xi_{s,l} \alpha_{s,l}.
\end{align*}
Thus $z$ is a linear combination of elements of $B$.
\epf

The proof of the next result follows by a straightforward computation.
\begin{lema}\label{yetter-struct-group} The $\ku G$-coaction of $S(\ku G, \ku_\psi F)$, given in \eqref{coact-coend-def}, is determined by
 $$\lambda(\alpha_{s,g}) = s^{-1}g s\ot \alpha_{s,g},$$
 for any $g\in F$, $s\in S$. The $\ku G$-action  on  $S(\ku G, \ku_\psi F)$, given in \eqref{H-mod-coend}, is determined by
 $$x\cdot \alpha_{s,g}=b(h^{-1},h^{-1}gh)\, \alpha_{r,h^{-1}gh},$$
 if $x=ft$ and $st^{-1}f^{-1}=hr$, where $f,h\in F$, $t,r\in S$.\qed
\end{lema}

For any subgroup $F\subset G $ and a normalized 2-cocycle $\psi\in Z^2(F,\ku^{\times})$, define  $C_\psi(G,F)$ as the subspace of $\Hom_\ku(\ku [S \times F], \ku)$ generated by functions $\phi:S \times F\to \ku$ such that
\begin{equation}\label{scalars-phi}
b(x, xs^{-1}gsx^{-1})\,\phi (s,g)=b(h^{-1}, h^{-1}gh)\,\phi(r,h^{-1}gh),
\end{equation}
for any $x\in G$  such that $sx^{-1}=hr$, $r\in S$, $h\in F$. Observe that if $F=G$, $\psi=1$ then $C_\psi(G,F)$ is the space of class functions on $G$.
\begin{prop}\label{class-funct-group} Let $F\subset G $ be a subgroup, and  $\psi\in Z^2(F,\ku^{\times})$ a normalized 2-cocycle. There exists a linear isomorphism  $\cf({}_{\ku_\psi F}\Mo)\simeq C_\psi(G,F)$.
\end{prop}
\pf As before, denote $S\subseteq G$ a set of representatives of elements of $F \backslash G$. Let be $\phi\in \cf({}_{\ku_\psi F}\Mo)$. This means that $\phi: S(\ku G, \ku_\psi F)\to S(\ku G, \ku G)$ is a morphism of $\ku G$-Yetter-Drinfeld modules. Elements of the basis described in Lemma \ref{basis} for $S(\ku G, \ku G)$ are of the form
$\alpha_{1,g}\in \ku^G\otk \ku G,$
for any $g\in G$. 

 Using Lemma \ref{yetter-struct-group}, since $\phi$ is a $\ku G$-comodule map, we observe that 
 $$\phi(\alpha_{s,g})= \phi_{s,g}\; \alpha_{1,s^{-1}gs},$$
for any $s\in S$, $g\in F$. Here $\phi_{s,g}\in \ku$. This implies that $\phi$ is determined by the scalars $\phi_{s,g}$. It remains to prove that these scalars satisfy \eqref{scalars-phi}. Take $x\in G$, and write it as $x=ft$, where $f\in F$, $t\in S$. Assume that $st^{-1}f^{-1}=hr$, where $h\in F$, $r\in S$. Since $\phi$ is a $\ku G$-module map, then
\begin{align*} \phi(x\cdot \alpha_{s,g})&=b(h^{-1},h^{-1}gh)\,\phi(\alpha_{r,h^{-1}gh})\\
&=b(h^{-1},h^{-1}gh)\,\phi_{r, h^{-1}gh}\; \alpha_{1,r^{-1}h^{-1}ghr}\\
&=x\cdot \phi(\alpha_{s,g})=\phi_{s,g}\; x\cdot \alpha_{1,s^{-1}gs}= \phi_{s,g}\;b(x, xs^{-1}gsx^{-1}) \alpha_{1,xs^{-1}gsx^{-1}}.
\end{align*}
This implies that $b(h^{-1},h^{-1}gh) \phi_{r, h^{-1}gh}=b(x, xs^{-1}gsx^{-1})\phi_{s,g}$.\epf

\subsection{Case $H=\ku^G$}\label{Section:case-dual-group} Let $F\subset G $ be a subgroup, and $\psi\in Z^2(F,\ku^{\times})$ be a normalized 2-cocycle. Let also $V$ be a simple $\ku_\psi F$-module. Recall the definition of the left $\ku^G$-comodule algebra $\kc(F,\psi,V)$ presented in Subsection \ref{Section:modgroup}. Again, let $S\subseteq G$ be a set of representatives of right cosets $F\backslash G$ such that $1\in S$. The following technical result will be needed later.

\begin{lema}\label{1dim-s} Let $f\in F$. The vector space consisting of $T\in \End(V)$ such that
\begin{equation}\label{commutation} U\circ T=T\circ (U\cdot f),
\end{equation} 
for any $U\in \End(V)$, is 1-dimensional.
\end{lema}
\pf Since the group $F$ is finite, the linear operator $f\cdot \, :V\to V$ is diagonalizable. Let $\{v_i\}_{i=1\dots n}$ be a basis of $V$ such that $f\cdot v_i=q_i\, v_i$, $q_i\in \ku^{\times}$ for any $i=1\dots n.$ Let $T\in \End(V)$  be a linear transformation such that it satisfies \eqref{commutation}. For any $j,k=1\dots n$ define $U_{j,k}:V\to V$ the operator $U_{j,k}(v_i)=\delta_{j,i} v_k$, for any $i=1\dots n$. Assume that 
$T(v_i)=\sum_{l} t_{i,l}\, v_l$. On one hand, for any $i=1\dots n$ we have that 
$$ (U_{j,k}\circ T)(v_i)=t_{i,j} v_k.$$
And, on the other hand we have
$$T \circ (U_{j,k}\cdot f)(v_i)=q_i q^{-1}_k \delta_{j,i} \sum_{l} t_{k,l}\, v_l.$$
Whence, equation \eqref{commutation} implies that if $i\neq j$, then $t_{i,j}=0$, and if $i=j$ then $q_i q^{-1}_k t_{k,k}=t_{i,i}$. This implies the Lemma.\epf

For any $f\in F$ denote by $T_f\in \End(V)$ the unique (up to scalar) non zero linear operator such that it fulfils condition \eqref{commutation} of Lemma \ref{1dim-s}.
For any $(f,s)\in F\times S$, denote $\alpha_{(f,s)}\in \ku G\otk \kc(F,\psi, V)$ by
$$\alpha_{(f,s)}= s^{-1} f s\ot \overline{ T_f\ot s}. $$
When $f=1$, we can choose $T_f=\Id_V$.
\begin{prop} The linearly independent set $\{\alpha_{(f,s)}:(f,s)\in F\times S\}$ is a basis for  $S(\ku^G, \kc(F,\psi,V)).$
\end{prop}
\pf It follows straightforward that for any $(f,s)\in F\times S$ the element $\alpha_{(f,s)}$ satisfies condition \eqref{sk-elements}. It follows from Lemma \ref{fpdim-char} (iv) that 
$$ \dim(S(\ku^G, \kc(F,\psi,V)))=\dim(\ku^G)=\mid G\mid.$$
Since the set $\{\alpha_{(f,s)}:(f,s)\in F\times S\}$ has cardinal equal to $\mid G\mid$, then it must be a basis.\epf

For any $(f,s)\in F\times S$ define $I(f,s)=\{(h,a)\in G\times G: aha^{-1}=s^{-1}fs\}$.
\begin{lema}\label{yetter-struct-dualgroup} The $\ku^G$-coaction of $S(\ku^G,  \kc(F,\psi,V))$, given in \eqref{coact-coend-def}, is determined by
$$\lambda(\alpha_{(f,s)})=\sum_{(a,h)\in I(f,s)} \delta_a\ot h\ot  \overline{ T_f\ot sa}$$
 for any $(f,s)\in F\times S$. The $\ku^G$-action  on  $S(\ku^G, \kc(F,\psi,V))$, given in \eqref{H-mod-coend}, is determined by
 $$ \delta_g\cdot \alpha_{(f,s)}=\begin{cases} 0 \quad \,\,\,\text{ if }\, g\neq s^{-1}fs\\
 \alpha_{(f,s)} \quad  \text{ if }\, g=s^{-1}fs.
 \end{cases}
 $$
\end{lema}\qed

\begin{prop} Let $F\subset G $ be a subgroup, and  $\psi\in Z^2(F,\ku^{\times})$ be a normalized 2-cocycle. Let $V$ be a simple $\ku_\psi F$-module. There exists a linear isomorphism  $\cf({}_{\kc(F,\psi,V)}\Mo)\simeq \ku^S$.
\end{prop}
\pf Recall that the comodule algebra representing the regular module category $\Rep(\ku^G)$ is the algebra $\kc(\{1\},1, \ku)$, see Remark \ref{regular-modcat}.

Let $\phi\in \cf({}_{\kc(F,\psi,V)}\Mo)$. Thus $\phi:S(\ku^G,  \kc(F,\psi,V))\to S(\ku^G,  \kc(\{1\},1, \ku)) $ is a $\ku^G$-module and $\ku^G$-comodule morphism. Elements in the basis of $S(\ku^G,  \kc(\{1\},1, \ku)) $ are $\alpha_{(1,g)}$, for any $g\in G$. Hence
\begin{align*}
\phi(\alpha_{(f,s)})&= \phi(\delta_{sfs^{-1}}\cdot\alpha_{(f,s)})=\delta_{sfs^{-1}}\cdot  \phi(\alpha_{(f,s)})\\
&=\begin{cases} 0 \quad \,\,\,\text{ if }\,  f\neq 1\\
\phi(\alpha_{(f,s)}) \quad  \text{ if }\, f=1.
 \end{cases}
\end{align*}
for any $f\in F$, $s\in S$. We deduce that $\phi$ is determined in values of $\alpha_{(1,s)}$ for any $s\in S$. Assume that 
$$\phi(\alpha_{(1,s)})=\sum_{g\in G} \phi_{s,g}\, \alpha_{(1,g)},$$
for certain $\phi_{s,g}\in \ku$. Since $I(1,s)=\{(1,a)\in G\times G\}$, and $\phi$ is a $\ku^G$-comodule map, we have that $\phi_{sa,g}=\phi_{s,ga^{-1}}$ for any $s\in S, g,a\in G$. Here we are abusing of the notation, since $sa$ denotes the element $t\in S$ that represents the class in which $sa$ belongs. Whence scalars $\phi_{s,g}$ are determined by a function $f: F\backslash G\to \ku$, as follows
$$ \phi_{s,g}=f(\overline{sg^{-1}}),$$
for any $s\in S, g\in G$.
\epf

\subsection{ The adjoint algebra for tensor categories $\ca(G,1,F,\psi)$}

Let $G$ be a finite group  and  $\omega\in Z^3(G,\ku^{\times})$ be a 3-cocycle. The category $\ca(G,\omega)$ stands for the category of  finite dimensional $G$-graded vector spaces, with associativity constraint defined by
$$a_{X,Y,Z}((x\ot y)\ot z)= \omega(g,h,f)\, x\ot (y\ot z),$$
for any $G$-graded vector space $X, Y, Z$, and any homogeneous elements $x\in X_g, y\in Y_h, z\in Z_f$.   Note that if $\omega=1$, then there is a monoidal equivalence $\ca(G,\omega)\simeq \Rep(\ku^G).$

If $F\subseteq G$ is a subgroup, and $\psi\in Z^2(F,\ku^{\times})$ is a 2-cocycle such that $d\psi\, \omega=1$, then the twisted group algebra $\ku_\psi F$ is an algebra in $\ca(G,\omega)$. The category $\ca(G,\omega,F,\psi)$ is the category ${}_{\ku_\psi F}\ca(G,\omega)_{\ku_\psi F}$ of $\ku_\psi F$-bimodules in $\ca(G,\omega)$. These categories are called \textit{group-theoretical fusion categories}.

We shall describe the adjoint algebra $\cha_{\Do}$, and the space of class functions $\cf(\Do)$ when $\Do=\ca(G,1,F,\psi)$. We shall keep the notation of Section \ref{Section:group-case}.

It follows from Corollary \ref{dual:caseHopf} that there is a monoidal equivalence
$$\ca(G,1,F,\psi)\simeq \Rep(G)^*_{{}_{\ku_\psi F}\Mo} $$
Using this equivalence and Lemma \ref{fpdim-char} (ii), it follows that $\cha_{\ca(G,1,F,\psi)}$ is isomorphic to $\theta_{\ku_\psi F}(\cha_{\ku_\psi F})$. Here $\cha_{\ku_\psi F}$ is the adjoint algebra corresponding to the $\Rep(G)$-module category ${}_{\ku_\psi F}\Mo$.

Using the explicit description of the functor $ \theta_{\ku_\psi F}$ given in Subsection \ref{Section:center-dual-Hopf}, we obtain that $\cha_{\ca(G,1,F,\psi)}$ is isomorphic to    $S(\ku G, \ku_\psi F)\otk  \ku_\psi F$. Recall that $S(\ku G, \ku_\psi F)$ has a basis consisting of elements  $\alpha_{s,l}\in\ku^G\otk \ku_\psi F$, $ s\in S, l\in F.$ The vector space $S(\ku G, \ku_\psi F)\otk  \ku_\psi F$ is an object in the category $\Zc({}_{\ku_\psi F}\ca(G,\omega)_{\ku_\psi F})$ as follows. The left $\ku G$-coaction  is given by
$$\lambda:S(\ku G, \ku_\psi F)\otk  \ku_\psi F\to \ku G\otk S(\ku G, \ku_\psi F)\otk  \ku_\psi F, $$
$$\lambda(\alpha_{s,f}\ot e_h)=s^{-1}fsh\ot \alpha_{s,f}\ot e_h,$$
for any $f,h\in F$, $s\in S$. The  $\ku_\psi F$-bimodule structure is given by
$$e_g\cdot  (\alpha_{s,f}\ot e_h)\cdot e_l= b(d^{-1}, d^{-1}fd) \psi(g,h)\psi(gh,l) \, \alpha_{r,d^{-1}fd}\ot e_{ghl},$$
if $sg^{-1}=df$, $g,f,l,d\in F$, $r,s\in S$.

The next result is a direct consequence of Lemma \ref{fpdim-char} and (the proof of) Proposition \ref{class-funct-group}.
\begin{lema} The space of class functions $\cf(\ca(G,1,F,\psi))$ is isomorphic to $C_1(G,F)$.\qed
\end{lema}

\end{document}